\documentclass[12pt]{amsart}
\usepackage{amsmath, amssymb, euscript, hyperref}

\newtheorem{thm}{Theorem}[section]

\newtheorem{prop}[thm]{Proposition}
\theoremstyle{definition}

\address{Azer Akhmedov, Department of Mathematics,
North Dakota State University,
Fargo, ND, 58108, USA}
\email{azer.akhmedov@ndsu.edu}

\begin{document}

 \begin{center} {\bf \Large Questions and Remarks on Discrete and Dense Subgroups of Diff(I)} \footnote{Future updates of this paper will be posted in the author's web page. Comments are welcome.}    \end{center}
  
  \bigskip

 \begin{center} {\bf Azer Akhmedov} \end{center}

\vspace{1cm}

 In recent decades, many remarkable papers  have appeared which are devoted to the study of finitely generated subgroups of $\mathrm{Diff} _{+}([0,1])$ (see {\bf [Be], [C], [FF1], [FF2], [FS], [G1], [G2], [N1], [N2], [Ts], [Y], [BW]} only for some of the most recent developments). In contrast, discrete subgroups of the group $\mathrm{Diff} _{+}([0,1])$ are much less studied. Very little is known in this area especially in comparison with the very rich theory of discrete subgroups of Lie groups which has started in the works of F.Klein and H.Poincar\'e in the 19th century, and has experienced enormous growth in the works of A.Selberg, A.Borel, G.Mostow, G.Margulis and many others in the 20th century. Many questions which are either very easy or have been studied a long time ago for (discrete) subgroups of Lie groups remain open in the context of the infinite-dimensional group $\mathrm{Diff} _{+}([0,1])$ and its relatives.
 
 \medskip
 
 In this survey article, we have collected questions about discrete and dense subgroups of $\mathrm{Diff} _{+}([0,1])$. Most of the questions are motivated by the problems discussed in {\bf [A1]-[A4]}. We have found these questions interesting and not trivial, although we believe that many of the questions we are asking are quite approachable.
 
 \medskip
 
 It is our pleasure to acknowledge the influences by many sources listed in the references. It is also a pleasure to thank to the participants of the {\em Dynamics Seminar} at North Dakota State University for their motivating questions.
 
 \medskip

 Throughout the paper, $I$ will denote the closed unit interval $[0,1]$. Let also
  
 $\mathrm{Homeo} _{+}(I)$ be the group of orientation preserving homeomorphisms of $I$;
 
 $\mathrm{Diff} _{+}(I)$ be the group of orientation preserving diffeomorphisms of $I$;
 
 $\mathrm{Diff} _{+}^{1+\epsilon }(I)$ be the group of orientation preserving diffeomorphisms of $I$ of regularity $C^{1+\epsilon }$ (a function $f:[0,1]\rightarrow \mathbb{R}$ belongs to class $C^{1+\epsilon}$ where $0<\epsilon <1$ if $f\in C^1$ and $|f'(x)-f'(y)|\leq M|x-y|^{\epsilon }$ for some $M > 0$ and for all $x,y\in I$).
 
 $\mathrm{Diff} _{+}^{\omega }(I)$ be the group of orientation preserving analytic diffeomorphisms of $I$.

  \medskip
  
  We study (discrete) subgroups of these groups by trying to treat them like Lie groups.  It turns out (justified by many strong results in the area) that while the groups $\mathrm{Homeo} _{+}(I)$ and $\mathrm{Diff} _{+}(I)$ are still quite {\em amorphous } (we borrow the term from {\bf [CK]}), the subgroups of $\mathrm{Diff} _{+}^{1+\epsilon }(I)$ start behaving more like subgroups of Lie groups. We limit ourselves here to mentioning only the following bright results:
  
  \medskip
  
  \begin{thm} a) (Plante-Thurston, {\bf [PT]}) Any nilpotent subgroup of $\mathrm{Diff} _{+}^2(I)$ is Abelian. 
  
  b) (Farb-Franks, {\bf [FF2]}) Any finitely generated torsion-free nilpotent group embeds in $\mathrm{Diff} _{+}(I)$.
  \end{thm}
  
  \medskip
  
  \begin{thm} (Navas, {\bf [N1]}) a) Any finitely generated subgroup of $\mathrm{Diff} _{+}^{1+\epsilon }(I)$ of sub-exponential growth is almost nilpotent.
  
  b) There exists a finitely generated group of intermediate growth which embeds in $\mathrm{Diff} _{+}(I)$. 
  \end{thm}
  
    \medskip
    
    Both Theorem 0.1. and Theorem 0.2. strongly indicate that while the group $\mathrm{Diff} _{+}(I)$ still exhibits a great flexibility and richness in terms of the combinatorial properties of its finitely generated subgroups, as soon as we demand a $C^{1+\epsilon }$ regularity, the behavior of the subgroups become much more tame. The fact that higher regularity implies some algebraic rigidity is a phenomenon known even from the older times dating back to at least Koppel {\bf [Kop]} and Thurston {\bf [Th]}.

  \medskip
  
  Instead of $I$, one can study homeomorphism/diffeomorphism groups of the other connected 1-manifolds: $\mathbb{R}, \ \mathbb{S}^1$ and $[0,1)$, but for the questions that we are interested in, the case of $I$ is the hardest and the most interesting. For some of the questions, we briefly sketch their analogue in the context of Lie groups, and provide references.

 \vspace{1cm}
 
 \section{\bf Free discrete subgroups}
 
 \medskip
   
  {\bf Question 1.} Does the free group $\mathbb{F}_2$ admit a faithful $C^1$-discrete representation in $\mathrm{Diff} _{+}^{1+\epsilon }(I)$?

  \bigskip
  
  It is known that $\mathbb{F}_2$ does admit a faithful $C^1$-discrete representation in $\mathrm{Diff} _{+}(I)$ ({\bf [A1]}). Question 1 (secretly important all along) comes into more light in connection with the recent very original approach of E.Shavgulidze to the amenability problem of subgroups of $\mathrm{Diff} _{+}(I)$ {\bf [Sh]}.   
  
  \medskip
  
  A $C_0$-discrete representation is clearly $C_1$-discrete (Mean Value Theorem). So one might ask if $\mathbb{F}_2$ admits a faithful $C_0$-discrete representation in $\mathrm{Diff} _{+}(I)$. The answer is negative and follows from Theorem A in {\bf [A2]}: {\em {\bf Theorem A.} Let $\Gamma $ be a subgroup of $\mathrm{Diff} _{+}(I)$, and $f, g\in [\Gamma ,\Gamma ]$ such that $f$ and $g$ generate a non-abelian free subsemigroup. Then $\Gamma $ is not $C_0$-discrete.}

  \medskip
  
  Containing a free semigroup is an extremely mild condition. Moreover, from practical point of view, it is usually very easy to prove that certain two elements in a group generate a free semigroup. (this may not be the case for free subgroups! For example, it is not easy to prove that the homeomorphisms $F(x) = x+1, G(x) = x^3$ of $\mathbb{R}$ indeed generate an isomorphic copy of $\mathbb{F}_2$ {\bf [Wh]}.)

  \medskip   
  
  Both $C_1$ and $C_0$ discreteness are very interesting in higher regularities. For $C_0$-discreteness, the following two theorems from {\bf [A2]} show that groups admitting $C_0$-discrete representation in higher regularities is a fairly small class.
  
  \medskip
  
  {\bf Theorem B.} {\em Let $\Gamma $ be a $C_0$-discrete subgroup of $\mathrm{Diff} _{+}^{1+\epsilon }[0,1]$. Then $\Gamma $ is solvable with solvability degree at most $k(\epsilon )$.}
     
     \medskip

     {\bf Theorem C.} {\em If $\Gamma $ is $C_0$-discrete subgroup of $\mathrm{Diff} _{+}^{2}[0,1]$ then $\Gamma $ is metaabelian.}

   \medskip
   
    Theorems B-C are obtained by combining Theorem A with the results of Navas {\bf [N1]}, Plante-Thurston {\bf [PT]}, and Szekeres {\bf [Sz]}. 
    
    \medskip
    
    {\bf Question 2.} What metaabelian subgroups of $\mathrm{Diff} _{+}^{2}[0,1]$ admit a $C_0$-discrete representation? 
 
  \medskip
  
  Question 2 can be asked about some specific and popular examples of metaabelian groups as well. One can probably hope for a classification result here.
  
  \medskip
  
  {\bf Question 3.} Let $f(x) = \frac{1}{2}x+\frac{1}{2}x^2$ and $g(x) = \frac{1}{3}x+\frac{2}{3}x^2$. The diffeomorphisms $f,g\in \mathrm{Diff} _{+}(I)$ generate a subgroup isomorphic to $\mathbb{F}_2$ (we will denote it by $\Gamma _{f,g}$). Is this subgroup $C_1$-discrete?
 
 \medskip
 
  We are inclined to think the answer should be negative but we are missing major steps of the proof. Discreteness and denseness can be thought of as opposite properties. The following question indicates how far we indeed are from the understanding of these concepts in the context of $\mathrm{Diff} _{+}(I)$:
  
  \medskip
  
  {\bf Question 4.} Is the subgroup $\Gamma _{f,g}$ $C_0$-dense in  $\mathrm{Diff} _{+}(I)$?
  
  \medskip

   Since the maps $f$ and $g$ are analytic, the following question becomes very interesting.
 
 \medskip
 
 {\bf Question 5.} Does $\mathbb{F}_2$ admit a faithful $C_1$-discrete representation in $\mathrm{Diff} _{+}^{\omega }(I)$?
 
 \medskip
 
  The main idea that we use for answering Question 3 is very useful also in here. But there are some general (plausible looking) claims that we have not succeeded in proving yet. As a side remark, notice that the subgroup generated by $f$ and $g$ is not $C_0$-discrete, by Theorem A.
 
 \medskip
 
 In contrast with Lie groups, very little study has been done of generic subgroups of $\mathrm{Diff} _{+}(I)$ (here,``generic" means that the group is generated by elements from an open dense subset). In the following two questions, by ``generic subgroup", we mean a subgroup generated by $k$ generic elements where $k > 1$.  
 
 \medskip
 
 {\bf Question 6.} Is a generic subgroup of $\mathrm{Diff} _{+}(I)$ $C_1$-discrete?  
 
 \medskip
 
 For a connected Lie group, random subgroups with sufficiently many generators turn out be dense (see {\bf [BG], [Wi2]} for the precise statements). Can we expect the same in $\mathrm{Diff} _{+}(I)$? More precisely,
 
 \medskip
 
 {\bf Question 7.} Is a generic subgroup of $\mathrm{Diff} _{+}(I)$ $C_0$-dense?
   
 \medskip
 
  It is interesting to compare Question 6 with Question 7; denseness and discreteness can be viewed as opposite properties. Our state of knowledge is so poor (this could well be the author's own ignorance) that we do not know what to expect from a generic subgroup. In fact, the situation is quite interesting even for Lie groups: let $G$ be a connected real Lie group, $k > 1$, and $D^k(G) = \{(g_1, \ldots , g_k)\in G^k \ | \ g_1, \ldots , g_k \ \mathrm{generate \ a \ discrete \ subgroup \ of} \ G\}$. It is proved in {\bf [Wi1]} that {\em both $D^k(G)$ and its compliment have positive measure when $G$ is not amenable. If $G$ is amenable, then the situation is also very interesting, namely, there exists $n$ such that $D^k(G)$ is of measure zero for $k > n$ and has complement of measure zero otherwise.}
  
  \medskip
     
    As emphasized in the previous section, the methods of {\bf [A2]} fail short in effectively understanding $C_1$-discreteness because although quite involved and subtle, the method essentially is based on the growth of certain subsets. However, just the understanding of $C_0$-discreteness has produced important results - we have applied Theorem B-C in solving an old problem extending the classical result of H\"older. 
    
    \medskip
    
    Let us recall a classical result (essentially due to H\"older,  cf.{\bf [N3]}) that if $\Gamma $ is a subgroup of Homeo$_{+}(\mathbb{R})$ such that every nontrivial element acts freely then $\Gamma $ is Abelian. A natural question to ask is what if every nontrivial element has at most $N$ fixed points where $N$ is a fixed natural number. In the case of $N=1$, we do have a complete answer to this question: it has been proved by Solodov (not published), Barbot {\bf [Ba]}, and Kovacevic {\bf [Kov]} that in this case the group is metaabelian, in fact, it is isomorphic to a subgroup of the affine group Aff$(\mathbb{R})$. (see {\bf [FF1]} for the history of this result, where yet another nice proof is presented). 
  
 \medskip
 
 Applying our results on $C_0$-discrete subgroups of Diff$(I)$, we have obtained the following results {\bf [A3]}.
  
  \medskip
  
  \begin{thm} Let $\Gamma $ be a subgroup of Diff $_{+}^{1+\epsilon }(I)$ such that every nontrivial element of $\Gamma $ has at most $N$ fixed points. Then $\Gamma $ is solvable.
  \end{thm}
  
  \medskip
  
  Assuming a higher regularity on the action we obtain a stronger result.
  
  \medskip
   
   \begin{thm} Let $\Gamma $ be a subgroup of Diff $_{+}^{2}(I)$ such that every nontrivial element of $\Gamma $ has at most $N$ fixed points. Then $\Gamma $ is metaabelian.
   \end{thm}
   
  \medskip
  
  An important tool in obtaining these results is provided by Theorems B-C from {\bf [A2]}. Theorem B (Theorem C) states that a non-solvable (non-metaabelian) subgroup of Diff $_{+}^{1+\epsilon }(I)$ (of Diff $_{+}^{2}(I)$) is non-discrete in $C_0$ metric. Existence of $C_0$-small elements in a group provides effective tools in tackling the problem. Such tools are absent for less regular actions, and the problem for Homeo$_{+}(I)$ (even for Diff$_{+}(I)$) still remains open.  
 
 \vspace{1cm}
 
 \section {\bf Margulis-Zassenhaus Lemma}
  
  \bigskip
  
  One of the major questions that we are interested in is the study of Margulis-Zassenhaus Lemma. This lemma (discovered by H.Zassenhaus in 1938, and later rediscovered by G.Margulis in 1968) states that in a connected Lie group $H$ there exists an open non-empty neighborhood $U$ of the identity such that any discrete subgroup generated by elements from $U$ is nilpotent (see {\bf [R]}). For example, if $H$ is a simple Lie group (such as $SL_2(\mathbb{R})$), and $\Gamma \leq H$ is a lattice, then $\Gamma $ cannot be generated by elements too close to the identity.  
     
     \bigskip
     
  In {\bf [A2]}, we prove weak versions of Margulis-Zassenhaus Lemma for the group $\mathrm{Diff} _{+}(I)$. It follows from the results of {\bf [A1]}, as remarked there, that the lemma does not hold either for $\mathrm{Diff} _{+}(I)$ or for $\mathrm{Homeo} _{+}(I)$, in $C_1$ and $C_0$ metrics respectively. 
     
     \medskip
     
   In the increased regularity the lemma still fails: given an arbitrary open neighborhood $U$ of the identity diffeomorphism in $\mathrm{Diff} _{+}(I)$, it is easy to find two $C^{\infty }$ ``bump functions" in $U$ which generate a discrete group isomorphic to $\mathbb{Z}\wr \mathbb{Z}$ which is not nilpotent; thus the lemma fails for $\mathrm{Diff} _{+}^{\infty }(I)$. 
   
    \medskip

  However, there is a strong evidence that the answer to the following question might indeed be positive.
 
 \medskip
 
  {\bf Question 8.} Does Margulis-Zassenhaus Lemma hold for $C_1$-discrete subgroups of $\mathrm{Diff} _{+}^{\omega }(I)$?
 
 \medskip
 
 Let us point out the following very interesting property of simple Lie groups (the proof of which uses Margulis-Zassenahus Lemma):
 
 \medskip
 
 {\em Stability of Generators:} We say a connected topological group $G$ has the stability of the generators property (SGP) if the following holds: if $g_1, \ldots , g_n$ generate a dense subgroup of $G$ then there exists an open non-empty neighborhood $U$ of identity such that for all $h_i\in g_iU, 1\leq i\leq n$, the subgroup generated by $h_1, \ldots , h_n$ is also dense in $G$. Using Margulis-Zassenhaus Lemma, it is indeed not difficult to prove this property for simple Lie groups.   
 
 \medskip
 
 {\bf Question 9.} Does $\mathrm{Diff} _{+}(I)$ have SGP in $C_1$-metric?   
    
 \medskip
 
  We can answer this question for $\mathrm{Homeo} _{+}(I)$ (hence also for $\mathrm{Diff} _{+}(I)$) in $C_0$-metric. The answer turns out to be negative. In fact, we prove {\bf [A4]} the following concrete statement which is interesting in itself (in the statement of the theorem, $||.||_0$ denotes the $C_0$ norm).
  
  \medskip
  
  \begin{prop} Let $\Gamma \leq \mathrm{Homeo} _{+}(I)$ be a finitely generated subgroup generated by homeomorphisms $g_1, \ldots , g_s$. Then for all $\epsilon > 0$, there exist homeomorphisms $h_1, \ldots , h_s$ such that $||h_i-g_i||_0 < \epsilon , 1\leq i\leq s$, and the group generated by $h_1, \ldots , h_s$ is $C_0$-strongly discrete, moreover, it is isomorphic to $\Gamma $. 
  \end{prop}
  
  \medskip
  
  Let us mention the definition of $C_0$-strong discreteness from {\bf [A1]}: {\em a subgroup $\Gamma $ is $C_0$-strongly discrete if there exists $\delta > 0$ and $x_0\in (0,1)$ such that $|g(x_0)-x_0| > \delta $ for all $g\in \Gamma \backslash \{1\}$.} 
  
  \medskip
  
  Similarly, one can define the notion of $C_1$-strong discreteness {\bf [A1]}: {\em a subgroup $\Gamma $ is $C_1$-strongly discrete if there exists $\delta > 0$ and $x_0\in (0,1)$ such that $|g'(x_0)-1| > \delta $ for all $g\in \Gamma \backslash \{1\}$.}
  
  \medskip
  
  Notice that a $C_0$-strongly ($C_1$-strongly) discrete subgroup is $C_0$-discrete ($C_1$-discrete). In light of Proposition 2.1, it is natural to ask the same question for the group $\mathrm{Diff} _{+}(I)$.
  
  \medskip
  
  {\bf Question 10.} Let $\Gamma \leq \mathrm{Diff} _{+}(I)$ be a finitely generated subgroup generated by diffeomorphisms $g_1, \ldots , g_s$. Then is there a $C_1$-strongly discrete subgroup (or just $C_1$-discrete subgroup) generated by diffeomorphisms from an artbitrarily small $C_1$-neighborhoods of $g_1, \ldots , g_s$? 
  
  \medskip 
  
 Proposition 2.1 implies that any finitely generated subgroup of $\mathrm{Homeo} _{+}(I)$ admits a $C_0$-strongly discrete faithful representation in $\mathrm{Homeo} _{+}(I)$. The following related question is borrowed from {\bf [A1]}:
 
 \medskip
 
 {\bf Question 11.} Is there a finitely generated group $\Gamma $ which admits a faithful representation in $\mathrm{Diff} _{+}(I)$ but does not admit a $C^1$-discrete faithful representation? 
  
   \medskip
   
   It is worth mentioning that, in a connected Lie group, the answer is always positive. For example, the group $\mathbb{Z}\wr \mathbb{Z}$ embeds in $GL(2,\mathbb{R})$ while it does not embed in discretely in any connected Lie group. 
     
 \medskip
    
  {\bf Question 12.} In regard to the Margulis-Zassenhaus Lemma, it is interesting to ask a reverse question, i.e. given an arbitrary open neighborhood $U$ of the identity in $G$, is it true that any finitely generated torsion free nilpotent group $\Gamma $ admits a faithful discrete representation in $\mathrm{Diff} _{+}(I)$ generated by elements from $U$? 
  
  \medskip
  
  In {\bf [FF2]}, it is proved that any such $\Gamma $ (i.e. any finitely generated torsion-free nilpotent group does admit a faithful representation into $\mathrm{Diff} _{+}(I)$ generated by diffeomorphisms from $U$. Also, it is proved in {\bf [N4]} that any finitely generated torsion-free nilpotent subgroup of $\mathrm{Diff} _{+}(I)$ indeed can be conjugated to a subgroup generated by elements from $U$. 
  
 \vspace{1cm}
 
 \section {\bf Dense subgroups of Diff(I)}
 
 \bigskip
 
  Discreteness and denseness can be viewed as opposite properties. Dense subgroups capture the (algebraic) content of the ambient group quite strongly while discrete subgroups may remain degenerate in the background of the ambient group unless we demand some uniformity. In the context of connected linear Lie groups, under slight technical conditions, a lattice (a discrete subgroup of finite covolume) turns out to be dense but in the (very weak) Zariski topology. Thus lattices (which also strongly capture the algebraic and geometric content of the ambient Lie group) can also be viewed as dense subgroups in a weaker topology.

  \medskip
  
  For Lie groups, (among the myriad of mostly well studied questions) it is often interesting to ask if a given finitely generated group $\Gamma $ admits a faithful discrete or dense representation. It is also interesting to ask if the same finitely generated group $\Gamma $ admits a faithful discrete representation and another faithful dense representation. (if $\Gamma = \mathbb{F}_2$, the answer to the latter question turns out be positive for any non-compact simple Lie group).

  \medskip
  
   For the group $\mathrm{Diff} _{+}(I)$, questions about discrete subgroups often time lead naturally to a question about dense subgroups (e.g. compare Questions 3 and 4, or see Question 6). In fact, since both discreteness and denseness are poorly understood, if a subgroup $\Gamma \leq \mathrm{Diff} _{+}(I)$ is non-discrete, one is often tempted to ask is $\Gamma $ perhaps dense?

 \medskip
 
   The following theorem indicates how far is the group $\mathrm{Diff} _{+}(I)$ from being solvable. (thus, $\mathrm{Diff} _{+}(I)$ cannot be possibly viewed as a relative of solvable Lie groups).
   
   \medskip
   
   \begin{thm} {\bf [A4]} A finitely generated $C_0$-dense subgroup $\Gamma $ of $\mathrm{Diff} _{+}(I)$ is not elementary-amenable.
   \end{thm}
   
   \medskip
   
   In the case of solvable groups, the proof of Theorem 3 is somewhat easier; we show that any finitely generated dense subgroup has infinite girth thus it cannot be solvable (see {\bf [A5]} for the definition of girth).
   
  \medskip 
   
   One approach to studying denseness is to consider an easier notion of dynamical transitivity. In fact, denseness can also be viewed as a dynamical $k$-transitivity, when $k = \infty $.

    \medskip
    
 We call the subgroup $\Gamma \leq \mathrm{Homeo} _{+}(I)$ on $I$ {\em $k$-transitive} if for all $z_1, \ldots , z_k \in (0,1)$, where $z_1 < z_2 < \ldots < z_k$ and for all open non-empty intervals $I_1, I_2, \ldots , I_k$ where $x < y$ if $x\in I_p, y\in I_q, 1\leq p < q\leq k$, there exists $\gamma \in \Gamma $ such that $\gamma (z_i)\in I_i, 1\leq i\leq k$. 
 
 \medskip
 
 Dynamical transitivity turns out to be quite an interesting notion. Besides serving as an intermediate step towards density, it is naturally interesting from dynamical point of view - transitivity of a group action may reveal information about the underlying algebraic structures of the group. But also, transitivity may play a role of a substitute for ergodicity (of discrete group actions on "homogenous spaces" of $\mathrm{Diff} _{+}(I)$; the latter cannot be defined easily in the absence of Haar measure in the infinite-dimensional group $\mathrm{Diff} _{+}(I)$).
 
 \medskip
 
 Dynamical transitivity is a key ingredient of the proofs of Theorem 1-2 from {\bf [A3]}. Proposition 1.5. of {\bf [A3]} proves 1-transitivity of irreducible non-solvable subgroups of $\mathrm{Diff} _{+}^{1+\epsilon }(I)$ under the condition that every non-trivial element has at most $N$ fixed points. The same method proves 1-transitivity for irreducible non-metaabelian subgroups of $\mathrm{Diff} _{+}^{\omega }(I)$ under no condition on the number of fixed points. 
 
 \medskip
 
 What about $k$-transitivity, for $k > 1$? Obtaining higher transitivity results may result in much bigger rewards.  We are still in the process of analyzing very elementary questions in this area. The following questions seem very interesting (in all these questions, we will assume that $\Gamma $ is an irreducible subgroup of $\mathrm{Diff} _{+}(I)$, i.e. it has no global fixed point):
 
 \medskip
 
 {\bf Question 13.} Is there a finitely generated $\Gamma \leq \mathrm{Diff} _{+}(I)$ which is $k$-transitive but not $(k+1)$-transitive?
 
 \medskip
 
 For $k = 1, 2$, one can easily find subgroups of Aff$(\mathbb{R})$ with the required property. It is interesting to study this question for $k\geq 3$.
 
  \medskip
  
  {\bf Question 14.} Let $k\geq 3$. For what values of $k$ does there exist a non-metaabelian subgroup $\Gamma \leq $ Diff $_{+}^{\omega }(I)$ which is $k$-transitive but not $(k+1)$-transitive?
  
  \medskip
  
  One can ask a similar question even for more specific subgroups:
  
  \medskip
  
  {\bf Question 15.} Let $\Gamma \cong \mathbb{F}_2 \leq $ Diff $_{+}^{\omega }(I)$. Is $\Gamma $ $k$-transitive for all $k$?
  
  \medskip

  As noted earlier, transitivity, as a dynamical property of a group action, may reveal algebraic properties of groups. It is especially intriguing to understand how transitive can a solvable or amenable group action be. We do not know the answer even to the following very simple question.
  
  \medskip
  
   {\bf Question 16.} Is there a finitely generated solvable $\Gamma \leq \mathrm{Diff} _{+}(I)$ which is $3$-transitive? Is there a $3$-transitive subgroup $\Gamma $ of $\mathrm{Diff} _{+}(I)$ which is not $C_0$-dense?
   
   \medskip

  In connection with Question 16, let us recall that finitely generated solvable groups cannot be $C_0$-dense. 
   
   \vspace{1cm}

   \section{\bf Further extensions of H\"older's Theorem}
   
  \medskip
  
  Theorems 1.1 and 1.2 pave a way to possibly further understanding of group actions on $I$ by homeomorphisms in terms of the sizes of the fixed point sets. 
  
  \medskip
  
  If $\Gamma \leq \mathrm{Homeo} _{+}(I)$ where all non-trivial elements have at mot $N$ fixed points then, for $N\leq 1$, one can indeed try to classify all such groups. For example, as noted earlier, we have the implications $$N=0\Rightarrow \Gamma \ \mathrm{is \ Abelian}$$ and $$N=1\Rightarrow \Gamma \ \mathrm{is \ isomorphic \ to \ a \ subgroup \ of} \ \mathrm{Aff}(\mathbb{R})$$
  
 \medskip
 
 {\bf Question 17.} Can one classify all subgroups of  $\mathrm{Homeo} _{+}(I)$  (or of $\mathrm{Diff} _{+}(I)$) such that every non-identity element has at most $N$ fixed points?
  
  \medskip

  If the number of fixed points of non-identity elements of $\Gamma $ is not uniformly bounded then it is interesting to know how fast this number grows. In fact, given any compact closed manifold $M$ with a diffeomorphism $T$ of $M$ with finitely many fixed points, it is natural to ask how fast the number of fixed points of $T^n$ grows as $n\to \infty $. If $a_n = |Fix(T^n)|, n = 1, 2, \ldots $, the study of the {\em the Artin-Mazur zeta function} $\zeta _f(z) = exp(\displaystyle \sum _{n=1}^{\infty }\frac{a_n}{n}z^n)$ has drawn considerable attention since 1960s ({\bf [AM], [Ka]}).  
  
  \medskip
  
  More generally, given a finitely generated group $\Gamma $ of homeomorphisms of $M$, it is natural to ask how fast the size of the fixed set $Fix(\gamma )$ grows when $n\to \infty $ where $\gamma \in \Gamma $ is an element of word length at most $n$ w.r.t. the fixed generating set.     
  
  \medskip
  
  For the next three questions, for simplicity, we restrict ourselves mostly to the analytic diffeomorphisms of $I$ but these questions make good sense also for any subgroup $\Gamma $ of $\mathrm{Diff} _{+}(I)$ where every non-identity element has finitely many fixed points. 
  
  \medskip
  
  Let $\Gamma \leq \mathrm{Diff} _{+}^{\omega }(I)$ be a finitely generated subgroup with a fixed finite generating set. Let also $$\omega _n(\Gamma ) = \mathrm{max} \{|Fix(\phi )| : \phi \in B_n(1), n = 1, 2, \ldots $$ where $B_n(1)$ denotes the ball of radius $n$ centered at the identity element $1\in \Gamma $.
  
  \medskip

  {\bf Question 18.} Can $\omega _n$ grow superexponentially? Can the growth of $\omega _n$ be sub-exponential but superpolynomial? 
  
  \medskip
  
  Since the actions of non-solvable subgroups of $\mathrm{Diff} _{+}^{1+\epsilon }(I)$ cannot admit a universal upper bound on the number of fixed points, we are wondering if this upper bound (i.e. the number $\omega _n$) may go to infinity very slowly. In particular,  
  
  \medskip
  
  {\bf Question 19.} Can one classify the groups $\Gamma $ for which $\omega _n$ grows linearly? Is there $\Gamma $ with a sub-linear growth of $\omega _n$? Can one possibly classify the groups $\Gamma $ for which $\omega _n$ grows polynomially?

  \medskip
  
  {\bf Question 20.} Compute the growth of $\omega _n$ for the group $\Gamma _{f,g}$.
  
  \medskip
  
  Since the universal bound on the number of fixed points implies solvability, it is interesting {\em \bf if a slow growth can be an indication of amenability}. In more general words, can the growth of $\omega _n$ detect amenability? Related to this, it is interesting to study the growth of $\omega _n$ for solvable subgroups of $\mathrm{Diff} _{+}(I)$. (For $\mathrm{Diff} _{+}^{\omega }(I)$, all solvable subgroups have been classified {\bf [BW]}).  
  
  \vspace{1cm}
  
  \section{\bf Lattices of Diff $_{+}^{1+\epsilon }(I)$}
  
  \medskip
  
  The study of discrete and dense subgroups naturally leads to the notion of a lattice. This notion makes better sense for the group $\mathrm{Diff} _{+}^{1+\epsilon }(I)$, although formally it can be defined for the group $\mathrm{Diff} _{+}(I)$ as well.
  
  \medskip
  
  {\em Definition. A subgroup $\Gamma \leq $Diff $_{+}^{1+\epsilon }(I)$ is called a lattice if it is $C_1$-discrete and $C_0$-dense.}
  
  \medskip
  
  A reader is invited to think of the similarities with lattices of linear connected Lie groups; their discreteness in the natural underlying Hausdorff topology and their density in the Zariski topology. Recall that, for a connected linear Lie group $G$, a lattice of $G$ is Zariski dense in it when a) $G$ is semisimple without compact factors {\em (Borel Density Theorem)},  or b) $G$ is simply connected and nilpotent {\bf [R]}. (The Zariski density result holds for a general connected linear Lie group but under slightly technical conditions.) Of course, a lattice is much more capturing than just being Zariski dense, but $C_0$ topology is also much richer than the Zariski topology. In {\bf [A4]}, we have studied some basic properties of lattices. We hope that, like in Lie groups, lattices of the group $\mathrm{Diff} _{+}^{1+\epsilon }(I)$ are very tightly related to it both algebraically and geometrically. 
  
  \medskip
  
  For example, in relation with the amenability problem, we expect that all lattices of $\mathrm{Diff} _{+}^{1+\epsilon }(I)$ are indeed non-amenable because the ambient group $\mathrm{Diff} _{+}^{1+\epsilon }(I)$ itself is non-amenable (we would like to caution that in the case of non-locally compact groups, the notion of amenability itself becomes somewhat more delicate). Notice that, by Theorem 3.1, an elementary amenable group cannot be a lattice of $\mathrm{Diff} _{+}^{1+\epsilon }(I)$. Notice also that, by {\bf [GS]}, R.Thompson's group $F$ has a faithful representation in $\mathrm{Diff} _{+}^{\infty  }(I)$ which is simultaneously $C_1$-discrete and $C_0$-dense. Hence non-amenability of $F$ should just follow from the fact that it is a lattice of a non-amenable group. At the moment, we do not have clear ideas to prove such a general and very strong conjecture but some basic facts understood on lattices of $\mathrm{Diff} _{+}^{1+\epsilon }(I)$ seem quite interesting and encouraging. For example, it is indeed true that lattices of $\mathrm{Diff} _{+}^{1+\epsilon }(I)$ cannot be elementary amenable.  

 \medskip
 
 {\bf Question 21.} Is every lattice of $\mathrm{Diff} _{+}^{1+\epsilon }(I)$   non-amenable?
 
 \medskip 
 
 The subgroups $\mathrm{Diff} _{+}^{\omega }(I)$ behave more rigidly and tamely. It is interesting therefore if a lattice of $\mathrm{Diff} _{+}^{1+\epsilon }(I)$ can sit totally inside a subgroup $\mathrm{Diff} _{+}^{\omega}(I)$. 
 
 \medskip
 
 {\bf Question 22.} Is there a subgroup of $\mathrm{Diff} _{+}^{\omega }(I)$ which is $C_1$-discrete and $C_0$-dense? 
 
 \medskip
 
   The group $\mathrm{Diff} _{+}(I)$, unlike $\mathrm{Diff} _{+}^{1+\epsilon }(I)$, is a huge source of amenable subgroups. But in $\mathrm{Diff} _{+}^{\omega }(I)$, we may witness that amenability implies a severe degeneracy.
   
   \medskip 
 
  {\bf Question 23.} Is it true that any subgroup of $\mathrm{Diff} _{+}^{\omega }(I)$ is either metaabelian or non-amenable?  
  
  \medskip
  
  In regard to let us emphasize that Tits Alternative for $\mathrm{Diff} _{+}^{\omega }(I)$ (asked by E.Ghys) is still not known. Many properties of groups can be passed from a Lie group to its lattice, and vice versa. One can ask Question 21 for properties other than amenability. Curiously enough, by {\em Thurston's Stability Theorem}, $\mathrm{Diff} _{+}(I)$ does not contain a subgroup with Kazhdan's property $(T)$. (a topological group is amenable and has property $(T)$ at the same time iff it is compact; so the two classes are essentially disjoint). Thus like the ambient group $\mathrm{Diff} _{+}(I)$, no lattice of it will have property $(T)$. Let us also recall a well known open question: does $\mathrm{Homeo} _{+}(I)$ have a subgroup with property $(T)$?

 \vspace{1cm}
 
 \section{\bf Divisible subgroups of Diff(I)}
 
 \bigskip
 
  The group $\mathrm{Diff} _{+}(I)$ is indeed a rich source of groups with interesting combinatorial properties. In questions about discreteness, one often encounters group elements which are infinitely divisible (i.e. an element has $n$-th root for any $n\geq 2$, or for infinitely many $n\geq 2$). It becomes interesting to understand if the $n$-th root converges to the identity (in $C^1$ metric or even in $C^0$-metric). 
  
  \medskip
  
  Let us recall that a group $G$ is called {\em divisible} if any element $g\in G$ has an $n$-th root for any natural $n$, i.e. there exists $f$ such that $f^n = g$. An element $h\in G$ is called {\em infinitely divisible} if for infinitely many $n\geq 2$ there exists $h_n\in G$ such that $h_n^n=h$.

  \medskip
  
  $\mathbb{Q}$ is a natural example of a divisible group although it is not easy to think a finitely generated divisible group (to our knowledge, the first such example has been constructed by V.Guba {\bf [Gu]}). 
  
  \medskip

  Despite exhibiting subgroups with interesting and even wild combinatorial properties, $\mathrm{Diff} _{+}(I)$ contains no finitely generated divisible subgroup. This immediately follows from {\em Thurston's Stability Theorem} which states that {\em any finitely generated subgroup of $\mathrm{Diff} _{+}(I)$ is indicable, i.e. it has a quotient isomorphic to $\mathbb{Z}$}. In {\bf [Be]}, it is shown that the group $G = \langle a, b, c \ | \ a^2 = b^3 = c^7 = abc\rangle$ embeds in  $\mathrm{Homeo} _{+}(I)$. On the other hand, $G$ is perfect. Thus Thurston's Stability Theorem fails for $\mathrm{Homeo} _{+}(I)$, and the following question becomes interesting.  
  
  \medskip
  
   {\bf Question 24.} Does $\mathrm{Homeo} _{+}(I)$ have a finitely generated divisible subgroup?    
  
  \medskip
  
  One is very tempted to conjecture that the answer to Question 24 should be negative, but it is not clear what tools to use. One might try to expand on the ideas of Thurston's proof but there are significant difficulties.
  
  \medskip
  
   {\bf Question 25.} Does $\mathrm{Diff} _{+}(I)$ have a finitely generated subgroup with an infinitely divisible non-trivial element $f$ and a sequence $f_{n_k}$ of $n_k$-th roots for an increasing sequence $n_k$ (i.e. $n_1 < n_2 < \ldots $ and $f_{n_k}^{n_k} = f, \forall k\in \mathbb{N}$) such that the sequence $f_{n_k}$ does not converge to identity in $C_0$-metric? 
  
  \medskip
  
  The answer to Question 25 is positive in $\mathrm{Homeo} _{+}(I)$, and probably negative in  $\mathrm{Diff} _{+}^{\omega }(I)$. For Lie groups, it is not difficult to see that there is no finitely generated divisible subgroup in a Lie group, moreover, the analogue of Question 25 has a negative answer. 
     
  \medskip
  
  A non-trivial homeomorphism has a square root in $\mathrm{Homeo} _{+}(I)$, in fact a whole continuum of them! (I have learned this from Matthew Brin). But a situation is quite different and in some instances very complicated in the groups  $\mathrm{Diff} _{+}(I)$ and $PL_{+}(I)$. Both the existence and uniqueness (of square roots, and more generally, of the $n$-th roots) questions are interesting. 
  
  \vspace{1cm}
  
 \section{\bf Subgroups of $PL_{+}(I)$} 
 
 In recent years, significant advances have been made in understanding (and in some instances even classifying) subgroups of PL$_{+}(I)$ from both algebraic and geometric points of view {\bf [BS], [Br1], [Br2], [Bl1], [Bl2]}. The representations of subgroups of PL$_{+}(I)$ in $\mathrm{Diff} _{+}(I)$ though are still very poorly understood.
 
 \medskip
 
 The following question is extremely intriguing.
 
 \medskip
 
 {\bf Question 26.} Does every finitely generated subgroup of PL$_{+}(I)$ admit a faithful representation in $\mathrm{Diff} _{+}(I)$?
 
 \medskip
 
 The answer is expected to be negative although we have very little idea at the moment about how can one approach to this problem. For a finitely generated group $\Gamma $, having a faithful embedding in $\mathrm{PL}_{+}(I)$ implies very strong structural rigidity from algebraic point of view. Despite this rigidity, $\mathrm{PL}_{+}(I)$ seems to be rich enough with finitely generated subgroups that not all of them have embeddings in $\mathrm{Diff} _{+}(I)$. The regularity of the embedding (if it exists) is also an interesting question. In fact it is useful to consider the following chain of subgroups.
 
 $$\mathrm{Homeo}_{+}(I)\supset \mathrm{Diff}_{+}(I) \supset \mathrm{Diff}_{+}^{1+\epsilon }(I)\supset \mathrm{Diff}_{+}^{2}(I)\supset \mathrm{Diff}_{+}^{\infty}(I)\supset \mathrm{Diff}_{+}^{\omega }(I)$$
 
 \ $\cup $
 
 $\mathrm{PL}_{+}(I)$
 
 \medskip
 
   For each of the inclusion in the above diagram, it is interesting to find finitely generated groups in the ambient group which do not have a faithful representation in the subgroup. Examples of finitely generated groups which embed in $\mathrm{Homeo}_{+}(I)$ but not in $ \mathrm{Diff}_{+}(I)$ are known {\bf [Be], [N1]}. On the other hand, Theorem 0.1 (combined results of Plante-Thurston and Farb-Franks) provide rich source of examples of finitely generated groups (namely, any f.g. torsion-free non-Abelian nilpotent group) embedding in $\mathrm{Diff}_{+}(I)$ but not in $\mathrm{Diff}_{+}^{2}(I)$. Also, combining Farb-Franks result with the result of A.Navas {\bf [N1]} one obtains examples of f.g. nilpotent groups embedding in $\mathrm{Diff}_{+}(I)$ but not in $\mathrm{Diff}_{+}^{1+\epsilon }(I)$. On the other hand, there are many examples of f.g. groups  which embed in $\mathrm{Diff}_{+}^{\infty}(I)$ but not in $\mathrm{Diff}_{+}^{\omega }(I)$. By the work of Ghys-Sergiescu {\bf [GS]}, R.Thompson's group $F$ also is a group like this. One can obtain a huge source of solvable groups with the required property using the classification result of Burslem and Wilkinson {\bf [BW]}.

   \medskip
   
   For the inclusions $G_1\supset G_2$ of the above diagram, it is also interesting to find a f.g. subgroup of $\mathrm{PL}_{+}(I)$ which embeds in the ambient group $G_1$ but does not embed in the subgroup $G_2$. The following question is interesting for any finitely generated group, not just for finitely generated subgroups of $\mathrm{PL}_{+}(I)$. 
 
 \medskip
 
 {\bf Question 27.} Is there a finitely generated group $\Gamma $ which emdes in $\mathrm{Diff}_{+}^{2}(I)$ but not in $\mathrm{Diff}_{+}^{3}(I)$?
 
 \medskip
 
 It is immediate to see that a non-Abelian subgroup of $\mathrm{PL}_{+}(I)$ is never $C_0$-discrete. But $C^1$-discreteness is not well understood at the moment. To be precise, let us define the following $C^1$-metric $d$ in PL$_{+}(I)$ : for any two maps $f, g\in \mathrm{PL}_{+}(I)$, let $$d(f, g) = \displaystyle \mathrm{sup}_{x\in (0,1)\backslash (S(f)\cup S(g))}|f'(x)-g'(x)|$$ where $S(h)$ denotes the (finite) set of singularities of a map $h\in \mathrm{PL}_{+}(I)$.
 
 \medskip
 
 Notice that, if the slope of $f\in \Gamma \leq \mathrm{PL}_{+}(I)$ at every point of $(0,1)\backslash S(f)$ belongs to a discrete multiplicative subgroup of $\mathbb{R}_{+}$ (e.g. R.Thompson's group $F$ in its standard representation in $\mathrm{PL}_{+}(I)$) then $\Gamma $ is $C_1$-discrete. Thus there are plenty of examples of interesting discrete subgroups of $\mathrm{PL}_{+}(I)$. However, unlike the case of $\mathrm{Diff} _{+}(I)$ it is apparently not easy to find  a non-discrete subgroup.       
 
 \medskip
 
 {\bf Question 28.} Is there a finitely generated non-discrete subgroup of   $\mathrm{PL}_{+}(I)$?

  \vspace{2cm}
  
  {\bf R e f e r e n c e s:}
  
  \bigskip
  
  [A1] Akhmedov, A. \ On free discrete subgroups of Diff(I). \ {\em Algebraic and Geometric Topology}, \ {\bf vol.4}, (2010) 2409-2418. 
  
  \medskip
  
  [A2] Akhmedov, A. \ A weak Zassenhaus Lemma for discrete subgroups of Diff(I). To appear {\em Algebraic and Geometric Topology.} \\ http://arxiv.org/pdf/1211.1086.pdf
  
  \medskip
  
  [A3] Akhmedov, A. \ Extension of H\"older's Theorem in $\mathrm{Diff} _{+}^{1+\epsilon }(I)$
  
  \medskip
  
  [A4] Akhmedov, A. \ On dense subgroups of Diff(I). Preprint.
  
  \medskip
  
  [A5] Akhmedov, A \ On the girth of finitely generated groups. \ {\em Journal of Algebra}, {\bf 268}, 2003, no.1, 198-208.
  
  \medskip
  
  [AM] Artin, M., Mazur, B. \ On periodic points \ {\em Annals of Mathematics}, {\bf 81}, (1965), 82-99.   
  
  \medskip
  
  [Ba] T.Barbot, \ Characterization des flots d'Anosov en dimension 3 par leurs feuilletages faibles, \ {\em Ergodic Theory and Dynamical Systems} {\bf 15} (1995), no.2, 247-270.
  
  \medskip
  
   [Be] Bergman, M. G.  \ Right orderable groups that are not locally indicable. \ {\em Pacific Journal of Mathematics.} {\bf 147} Number 2 (1991), 243-248.
		
		\medskip
		
		[BG] Breuillard, E., Gelander, T. \ On dense free subgroups of Lie groups. \ {\em Journal of Algebra}, {\bf vol.261}, no.2, 448-467.
		
		\medskip
		
		[Bl1] Bleak, C. \ A geometric classification of some solvable groups of homeomorphisms. {\em Journal of London Mathematical Society} (2), {\bf 78}, (2008) no. 2, 352-372.  
		
		\medskip
		
		[Bl2] Bleak, C. An algebraic classification of some solvable groups of homeomorphisms. \ {\em Journal of Algebra.} {\bf 319}, 4, p. 1368-1397
		
		\medskip
		
		[Br1] Brin, M. \ The ubiquity of Thompson's group F in groups of piecewise linear
homeomorphisms of the unit interval, \ {\em Journal of London Math. Soc.} (2), {\bf 60}  (1999), no. 2, 449–460.

		\medskip
		
		[Br2] Brin, M. \ Elementary amenable subgroups of R. Thompson's group F. \ {\em International Journal of Algebra and Computation,} {\bf 15} (2005) no 4, 619-642.
		
		\medskip
		
		[BS] Brin, M., Squier, C. \ Groups of piecewise linear homeomorphisms of the real line.  \ {\em Inventiones Mathematicae} {\bf 79} (1985), no.3, 485-498.
		
		\medskip
		
		[BW] Burslem, L. Wilkinson, A. \  Global rigidity of solvable group actions on $\mathbb{S}^1$. \ {\em Geometry and Topology}, {\bf 8},  (2004), 877--924
		
		\medskip
		 
     [C] Calegari, D.  \ Nonsmoothable, locally indicable group  
actions on the interval. \ {\em Algebraic and Geometric Topology} {\bf 8} (2008) no. 1, 609-613     

  \medskip
	
	   [CK] Cohen. M., Kallman, R. \ Openly Haar null sets and conjugacy in Polish groups. \  {\em Preprint.}
	   
	 \medskip
	    
	   [FF1] Farb, B., Franks, J. \ Groups of homeomorphisms of one-manifolds, II: Extension of Holder's Theorem. {\em Transactions of the  American Mathematical Society.} {\bf 355}, (2003) , no.11, 4385-4396.
		
		\medskip
		
		[FF2] Farb, B., Franks, J. \ Groups of homeomorphisms of one-manifolds, III: Nilpotent subgroups. \ {\em Ergodic Theory and Dynamical Sysytems.} {\bf 23}, (2003), no.5, 1467-1484. 
		
		\medskip
		
		[FS] Farb, B., Shalen, P. \ Groups of real-analytic diffeomorphisms of the circle. \ {{\em Ergodic Theory and Dynamical Sysytems.} {\bf 22}, (2002), 835-844.
		
		\medskip
		
	   [G1] Ghys, E. \ Groups acting on the circle. \ {\em Enseign. Math. (2)} 47(3-4):329-407, 2001.
		
		\medskip
		
		 [G2] Ghys, E. \ Sur les groupes engendrés par des difféomorphismes proches de l'identité.  Bol. Soc. Brasil. Mat. (N.S.) 24 (1993), no. 2, 137-178.
		 
		 \medskip

		[GS] Ghys, E., Sergiescu, R. Sur un groupe remarquable de diffeomorphismes du cercle, Comment.Math.Helvetici 62 (1987) 185-239.

	 \medskip
	 
	 [Gu] Guba, V. A finitely generated complete group. \ Izv. AN SSSR, 1986, {\bf 50}, no. 5, 883-924.
	 	 
	 \medskip
	 
	 [Ka] Kaloshin, V. \ An extension of the Artin-Mazur Theorem. \ {\em Annals of Mathematics}, {\bf 150} (1999), 729-741.
	 
	 \medskip
	 
	 [Kop] Kopell, N. Commuting Diffeomorphisms. \ In Global Analysis, Pros. Sympos. Pure Math XIV Berkeley, California, (1968) 165-184.
	  
	 \medskip
	 
	 [Kov] Kovacevic, N. \ M\"obius-like groups of homeomorphisms of the circle. \ {\em Trans. Amer. Math. Soc.} {\bf 351} (1999), no.12, 4791-4822.
	 
	 \medskip
	 	
      [N1] Navas, A. \ Growth of groups and diffeomorphisms of the interval. \ {\em GAFA.} {\bf 18}, (2008), 988-1028.
					
		\medskip
		
		 [N2] Navas, A. \ A finitely generated locally indicable group with no faithful action by $C^1$ diffeomorphisms of the interval. {\em Geometry and Topology} {\bf 14} (2010) 573-584.
			  
   \medskip
   
   [N3] Navas, A. Groups of Circle Diffemorphisms. \ Chicago Lectures in Mathematics, 2011. {\em http://arxiv.org/pdf/0607481} 
   
   \medskip
  
   [N4] Navas, A. Sur les rapprochements par conjugasion en dimension 1 et classe $C^1$. \ {\em http://arxiv.org/pdf/1208.4815}.
  
  \medskip

 [PT] Plante, J., Thurston, W. \ Polynomial growth in holonomy groups of foliations. \ {\em Comment. Math. Helv.} \ {\bf 51} (1976), 567-584. 
  
  \medskip

 [R] Raghunathan, M.S. {\em Discrete subgroups of semi-simple Lie groups}. \ Springer-Verlag, New York 1972. Ergebnisse der Mathematik und ihrer Grenzgebiete, Band 68. 
  
  \medskip
  
  [Sh] Shavgulidze, E.T. \ About amenability of subgroups of the group of diffeomorphisms of the interval. \ http://arxiv.org/pdf/0906.0107.pdf
  
  \medskip
  
  [Sz] Szekeres, G. \ Regular iteration of real and complex functions. {\em Acta Math.} \ {\bf 100} (1958), 203-258     
  
  \medskip
  
  [Th]  Thurston, W. \ A generalization of Reeb stability theorem. \ {\em Topology} {\bf 13}, (1974) 347-352. 
     
     \medskip
		
	[Ts] Tsuboi, T. \ Homological and dynamical study on certain groups of Lipschitz homeomorphisms of the circle. \ {\em Journal of Mathematical Society of Japan.} {\bf 47} (1995) 1-30.
		
    \medskip
		
	[Y] Yoccoz, J.C. \ Centralisateurs et conjugasion diff\'erentiable des diff\'eomorphismes du cercle. Petits diviseurs en dimension 1. {\em Ast\'erisque} {\bf 231} (1995), 89-242.
  
  \medskip
  
 [Wh] White, S. $X\rightarrow X+1$ and $X\rightarrow X^p$  generate a free subgroup. Journal of Algebra, {\bf 118}, (1988), 408-422.

 \medskip
 
 [Wi1] Winkelmann, J. \ Generic subgroups of Lie groups. \ {\em Topology}, {\bf 41}, (2002), 163-181.
 
 \medskip
 
 [Wi2] Winkelmann, J. \ Dense random finitely generated subgroups of Lie groups. \  http://arxiv.org/abs/math/0309129

\end{document}